\documentclass[preprint,10pt]{elsarticle}
\usepackage{epsfig,array,graphicx}
\usepackage{amsmath,amssymb}
\usepackage{rotating,rotfloat}
\usepackage{multirow}
\usepackage{subfigure}
\begin{document}

\begin{frontmatter}

\title{Solving systems of transcendental equations involving the Heun functions.}
\author[a1,a2]{Plamen P. Fiziev\corref{cor1}}
\ead{fiziev@phys.uni-sofia.bg, fiziev@theor.jinr.ru}
\author[a1]{Denitsa R. Staicova}
\ead{dstaicova@phys.uni-sofia.bg}
\cortext[cor1]{Corresponding author}
\address[a1]{Department of Theoretical Physics, Sofia University ``St. Kliment Ohridski", 5
James Bourchier Blvd., 1164 Sofia, Bulgaria}
\address[a2]{BLTF, JINR, Dubna, 141980 Moscow Region, Russia}

\begin{abstract}
The Heun functions have wide application in modern physics and are expected to succeed the hypergeometrical functions in the physical problems of the 21st century. The numerical work with those functions, however, is complicated and requires filling the gaps in the theory of the Heun functions and also, creating new algorithms able to work with them efficiently. 

We propose a new algorithm for solving a system of two nonlinear transcendental equations with two complex variables based on the M\"uller algorithm. The new algorithm is particularly useful in systems featuring the Heun functions and for them, the new algorithm gives distinctly better results than Newton's and Broyden's methods.

As an example for its application in physics, the new algorithm was used to find the quasi-normal modes (QNM) of Schwarzschild black hole described by the Regge-Wheeler equation. The numerical results obtained by our method are compared with the already published QNM frequencies and are found to coincide to a great extent with them. Also discussed are the QNM of the Kerr black hole, described by the Teukolsky Master equation.

\end{abstract}

\maketitle

\end{frontmatter}

\section{Overview}
The Heun functions appear with increasing frequency in modern physics. For example, they arise in the Schr\"odinger equation with anharmonic potential, in water molecule, in the Stark effect,  in different quantum phenomena related with repulsion and attraction of levels, in the  theory of lunar motion, in gravitational physics of scalar, spinor, electromagnetic and gravitational waves in Schwarzschild and Kerr metric, in crystalline materials, in three-dimensional waves in atmosphere, in Bethe anzatz systems, in Collogero-Moser-Sutherland systems, e.t.c., just to mention a few. Because of the wide range of their applications (\cite{DLMF,heun3_})-- from quantum mechanics to astrophysics, from lattice systems to economics  -- they can be considered as the $21^{\text{st}}$ century successors of the hypergeometric functions encountered in some simple physical problems of $20^{\text{th}}$ century. 

It is not hard to explain this situation. In natural sciences, in particular, in physics we usually study the different phenomena starting from some equilibrium state. Then we study small deviations from it in linear approximation, and at the end, going far away from the equilibrium we are forced to take into account nonlinear phenomena. It is well known that to describe the wave processes (like those in quantum mechanics), related with some linear phenomenon in classical physics (like classical mechanics), we have to use hypergeometric functions. Therefore these functions were well studied in $19^{\text{th}}$ and $20^{\text{th}}$centuries and today one can find the corresponding codes in all good computer packages. According to the theorem by Slavyanov \cite{heun3_}, if we study nonlinear classical phenomena, described by elliptic functions, or even by the solutions of any of Painlev\'e type equations, the corresponding wave problems can be solved exactly in terms of the Heun functions. Since the Painlev\'e equations can be considered as Hamilton ones for a very large class of nonlinear classical problems, one can expect a fast increase in the applications of the Heun functions in physics and other natural sciences of $21^{\text{th}}$ century.       

Their mathematical complexity, however, makes working with them a significant challenge both analytically and numerically. The Heun functions are unique local Frobenius solutions of a second-order linear ordinary differential equation of the Fuchsian type \cite{heun3_,heun,heun1_,heun2_} which in the general case have 4 regular singular points. Two or more of those regular singularities can coalesce into an irregular singularity leading differential equations of the confluent type and their solutions: confluent Heun function, biconfluent Heun function, double confluent Heun function and triconfluent Heun function. The Heun functions generalize the hypergeometric function (and also include the Lam\'e function, Mathieu function and the spheroidal wave functions \cite{heun3_,heun2_}) and some of their uses can be found in \cite{heun3_} and also in the more recent \cite{heun1}. Clearly, the Heun functions will be encountered more and more in modern physics, hence, there is a need of better understanding of those functions and new, more adequate algorithms working with them.

Despite the growing number of articles which use equations from the Heun type and their solutions, the theory of those functions is far from complete. There are some analytical works on the Heun functions, but they were largely neglected until recently. Some recent progress can be found in \cite{Fiziev4}, but as a whole there are many gaps in our knowledge of those functions. Particularly, the connection problem for the Heun functions is not solved -- one cannot connect two local solutions at different singular points using known constant coefficients \cite{heun3_}.   Another example of a serious gap in the general theory of the Heun functions in general is the absence of integral representations analogous to the one for hypergeometric functions. According to Whittaker's hypothesis, the Heun's functions are the simplest class of special functions for which no representations in form of contour integrals of elementary functions exists. In some particular cases are known integral representations in form of double integrals of elementary functions, but the general case is an open problem.

Numerically, the only software package currently able to work with the Heun functions is \textsc{maple}. Alternative ways for evaluations of those functions do not exist (to the best of our knowledge) and there are no known projects aiming to change this situation, an admittedly immense task by itself. This means that the use of the Heun functions is limited to the routines hidden in the kernel of \textsc{maple}, which the user cannot change or improve -- a situation that makes understanding the numerical problems or avoiding them adequately very difficult. On the positive side, those routines were found by the team to work well enough in many cases (see, for example, the match between theory and numerical results in  \cite{spectra}, as well as the other applications of those functions in \cite{arxiv3,arxiv_Kerr}). Yet, there are some peculiarities -- there are values of the parameters where the routines break down leading to infinities or to numerical errors. The situation with the derivatives of the Heun functions in \textsc{maple} is even worst -- for some values they simply do not work, for example outside the domain $|z|<1$, where their precision is much lower than that of the Heun function itself. Also, in some cases there are no convenient power-series representations and then the Heun functions are evaluated in \textsc{maple} using numerical integration. Therefore the procedure goes slowly in the complex domain (compared to the hypergeometric function) which means that the convergence of the root-finding algorithm is essential when one solves equations including Heun's functions. 

Despite all the numerical set-backs, the Heun functions offer many opportunities to modern physics. They occur in the problem of quasi-normal modes (QNM) of rotating and nonrotating black holes, which is to some extent the gravity analogue of the problem of the hydrogen atom. Finding the QNMs is critical to understanding observational data from gravitational wave detectors and proving or refuting the black holes existence (\cite{QNM,QNM2} and also \cite{spectra,arxiv3,arxiv_Kerr}). In this case, one has to solve a two-dimensional connected spectral problem with two complex equations in each of which one encounters the confluent Heun functions. The analytical theory of the confluent Heun function is more developed than that of the other types of Heun functions, but still many unknowns remain. Again, the evaluation of the derivative of the confluent Heun function is problematic in \textsc{maple}, which makes Newton's root-finding algorithm (\cite{Muller,numerical}) unusable. Broyden's algorithm (\cite{Broyden}) generally works, but it is slowly convergent even close to a root (see \cite{arxiv}). It is clear that we need a novel algorithm, that will offer quicker convergence than Broyden's algorithm, but without relying on derivatives.

To solve this problem, in the case of a system of two equations in two variables, our team developed a two-dimensional generalization of the M\"uller algorithm.  The one-dimensional M\"uller algorithm (\cite{Muller2}) is a quadratic generalization of the secant method, that works well in the case of a complex function of one variable. It has very good convergence for a large class of functions ($\sim 1.84$) and it is very efficient when the starting point (the initial guess) is close to a root (for applications see \cite{Fiziev1} and \cite{spectra}). It is also well convergent when working with special transcendental functions such as the Heun functions. Most importantly, this algorithm does not use derivatives, which is key for our work.  

The two-dimensional M\"uller algorithm seems to inherit some of the advantages of its one-dimensional counterpart like good convergence and usability on large class of functions as our tests show. The new algorithm was used to solve the QNM problem in the case of a Schwarzschild black hole and it proved to work without significant deviations from the results published by Andersson (\cite{Q_N_M}) and Fiziev (\cite{Fiziev1}).  Also, preliminary results for the QNM of the Kerr black hole are discussed and for them we also obtain a very good coincidence with published results \cite{special31}.

The article is organized as follows: Section 2 reviews the one-dimensional M\"uller algorithm and it introduces in detail its two-dimensional generalization, in Section 3 we give two examples of the application of the new method in systems featuring the confluent Heun functions -- the QNMs of rotating and nonrotating BH, and the results in both cases are analyzed in terms of precision and convergence. In Section 4 we summarize our results.

\section{The M\"uller algorithm}
\subsection{One-dimensional M\"uller's algorithm}
The one-dimensional M\"uller algorithm (\cite{Muller,Muller2}) is iterative method which at each step evaluates the function at three points, builds the parabola crossing those points and finds the two points where that parabola crosses the x-axis.  The next iteration is the the point farthest from the initial point.

Explicitly, the one-dimensional M\"uller algorithm can be defined as the map $ \mu: \quad \mu(x^{0},F(x))\rightarrow x^{P}$, which obtains the final point $x^P$ in $P$ iterations by calculating for every three points $x_{{j-2}}$, $x_{{j-1}}$, $x_{{j}}$ (with the corresponding values of the function $f(x)\rightarrow f_{{j-2}}$, $f_{{j-1}}$, $f_{{j}}$), the next iteration $x_{{j+1}}$ as:
\begin{align*}
&x_{{j+1}}=x_j\!-\left( x_{{j}}\!-\!x_{{j-1}} \right)\frac{2C}{{{\it max(D_1,D_2)}}},\, \text{where}\\
&A=f_j\,q\!-\!q(1\!+\!q)f_{j-1}+q^2 f_{j-2},\, B=\left( 2\,q\!+\!1 \right)f_{{j}} \! -\! \left( 1\!+\!q \right) ^{2}f_{{j-1}}\!+\!{q}^{2}f_{{j-2}},\\
&C= \left( 1\!+\!q \right) f_{{j}},\, {\it 	D_{1,2}}=B\pm\sqrt {{B}^{2}\!-\! 4\,AC}\,\text{and} \,\,q={\frac {x_{{j}}-x_{{j-1}}}{x_{{j-1}}-x_{{j-2}}}}. 
\end{align*}

\noindent The iterations continue until $\mid\!x_{P}-x_{P-1}\!\mid< 10^{-d},$ where $d$ is the number of digits of precision we require. This exit-condition works independently of the actual numerical zero in use, which may vary for the confluent Heun function and thus it is the most appropriate for our numerical work.

\subsection{Two-dimensional M\"uller's algorithm}
The two-dimensional M\"uller method comes as a natural extension of the one-dimensional M\"uller method.

For two complex-valued functions $F_1(x,y)$ and $F_2(x,y)$ we want to find such pairs of complex numbers $(x_{{}_I},y_{{}_I})$ which are solutions of the system:
\begin{equation}
\begin{cases}
F_1(x_{{}_I},y_{{}_I})=0\\
F_2(x_{{}_I},y_{{}_I})=0
\end{cases}
\label{sys1}
\end{equation}
\noindent where $I=1,\ldots$ numbers the solution in use. From now on, we will omit the index $I$, considering that we work with one arbitrary particular solution. Finding all the solutions of a system is beyond the scope of this article.

Consider the functions $F_1(x,y),\, F_2(x,y)$ as two-dimensional {\em complex} surfaces $z=F_1(x,y)$ and $z=F_2(x,y)$ in a three-dimensional space of the {\em complex} variables $\{x,y,z\}$\footnote{Equivalently, we can consider four {\em real} surfaces in five-dimensional {\em real} hyperspace, which are defined by four {\em real} functions of four {\em real} variables}. Normally, to solve the system, one expresses the relation $y(x)$ from one of the equations, then by substituting it in the other equation, one solves it for $x$ and from $y(x)$ one finds $y$. In the general case, however, this is not possible. The idea of our code is to approximately follow that procedure by finding an approximate linear relation $y(x)$ between the two variables and then using it to find the root of function of one variable trough the one-dimensional M\"uller algorithm.

\begin{figure}
\vspace{-0cm}
\hspace{1.5cm}
\includegraphics[width=300px,height=200px]{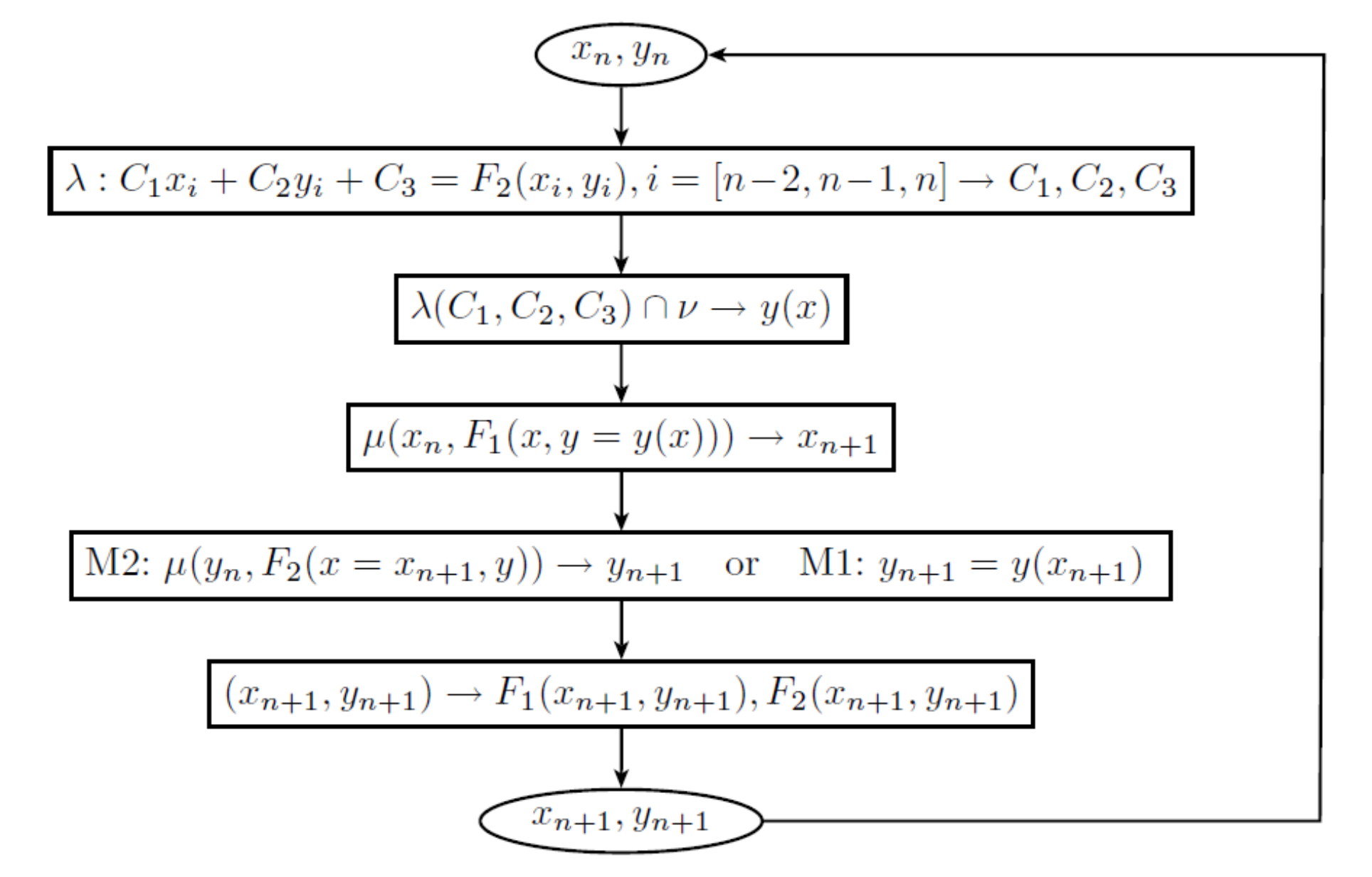}
\caption{A block scheme of the two-dimensional M\"uller algorithm. $\lambda(C_1, C_2, C_3)$ is the plane with equation $z=C_1x + C_2y + C_3$ that crosses trough the 3 pairs of points $(x_{i},y_{i})$ and the function $F_2$ evaluated in them. The plane $\nu$ is defined by the equation $z=0$. The one-dimensional M\"uller algorithm, $\mu(t^{in},F(t)) \rightarrow t^{fin}$, is applied on the function of {\em one variable} $F(t)$ with starting point $t^{in}$ and final point $t^{fin}$.}
\label{block}
\end{figure}

To find the linear relation $y(x)$, at each iteration we form the plane passing trough three points of one of the functions and then the equation of the line of intersection between that plane and the plane $z=0$ is used as the approximate relation $y(x)$. This basically means that the so found $y(x)$ is an approximate solution of one of the equations which ideally should be near the real solution in the $z=0$ plane. Substituting this relation in the other function, we run the one-dimensional M\"uller algorithm on it to fix the value of one of the variables, say $x$. Using the value of $x$ in the first function, we again run the one-dimensional M\"uller algorithm on it to fix the value of the other variable -- $y$. Alternatively one can substitute the value of $x$ directly in $y(x)$ to obtain $y$. This ends one iteration of the algorithm.  The process repeats until one of the exit-conditions discussed below has been reached. The block-scheme of the algorithm can be seen on Fig. (\ref{block}).

Explicitly, the code starts by evaluating the two functions $F_{1, 2}(x_i, y_i)$ in three starting pairs of points ($i=1,2,3$) that ideally should be near one of the roots of the system. In our case, those three \emph{initial} pairs are obtained from one starting pair to which we add and subtract certain small complex number. This artificial choice is done only in the first iteration ($n=3$), afterwards we use the output of the last three iterations to form $(x_{n-2},y_{n-2}), (x_{n-1},y_{n-1})$, $(x_n,y_n)$ and the respective $F_{1,2}(x,y)$. Thus on every iteration after $n=3$ the actual complex functions $F_{1, 2}(x_n, y_n)$ are evaluated only once outside of the one-dimensional M\"uler subroutines.

Next we construct the plane passing trough those three points for one of the functions, say $F_2$ by solving the linear system:
\begin{align*}
&C_1x_{n-2}+C_2 y_{n-2}+C_3=F_2(x_{n-2},y_{n-2})\\
&C_1x_{n-1}+C_2 y_{n-1}+C_3=F_2(x_{n-1},y_{n-1})\\
&C_1x_{n}+C_2 y_{n}+C_3=F_2(x_{n},y_{n}).
\end{align*}
\noindent From it one obtains the coefficients $C_1, C_2, C_3$ of the plane $z\!=\!C_1x+C_2y+C_3$.

This plane is intersected with the plane $z=0$ and the equation of the line between those two planes is the approximate relation $y(x)$ (i.e. $y=(-C_1x-C_3)/C_2$) of the two variables.

We substitute that relation in the first function $F_1(x,y) \rightarrow F_1(x,y(x))$ and we start the one-dimensional M\"uller on that ``linearized`` function of only one variable, $x$. After some pre-determined maximal number of iterations, the exiting point is chosen for $x_{n+1}$ ( $\mu(x_n,F_1(x,y=y(x))) \to x_{n+1}$) .

Then, there are two possibilities.

Algorithm M1: one could use directly the relation $y(x=x_{n+1})$ to find $y=y_{n+1}$. Or,

Algorithm M2: One can substitute $x=x_{n+1}$ in the other function $F_2(x,y) \to F_2(x=x_{n+1}, y)$ in order to find $y_{n+1}$ using again the one-dimensional M\"uller algorithm ($\mu(y_n,F_2(x_{n+1},y)) \to y_{n+1}$).

Our numerical experiments showed that both approaches lead to convergent procedure.

After $(x_{n+1},y_{n+1})$ are fixed, the two functions $F_{1, 2}(x_{n+1}, y_{n+1})$ are evaluated and if the new points are not roots, the iterations continue.

The exit-strategy in the two-dimensional M\"uller algorithm is as follows:
\begin{enumerate}
 \item To avoid hanging of the algorithm or its deviation from the actual root of the system, one fixes the maximal number of iterations for the one-dimensional M\"uller subroutine, $P$. From our experience, $P\!=\!4\!-\!10$ gives best convergence.
 \item The precision-condition ($\mid\! x_{j}\!-\!x_{j-1}\!\mid<10^{-d}$) remains in force for the one-dimensional M\"uller algorithm. Usually the subroutine exits, because of $j>P$ during the first few iterations of the two-dimensional M\"uller algorithm, until it gets closer to the roots and exits with $\mid\! x_{j}\!-\!x_{j-1}\!\mid<10^{-d}$.
 \item The primary exit-condition of the two-dimensional M\"uller is fulfilled when for two consecutive pairs $(x_{n}, y_{n})$: $\mid x_{n} - x_{n-1}\mid < 10 ^{-d}, \mid y_{n} - y_{n-1}\mid < 10 ^{-d}$.  In this case, if the values of functions $F_1(x,y), F_2(x,y)$ at those points are sufficiently small, the algorithm exits with a root.
 \item To keep the two-dimensional M\"uller algorithm from hanging, a maximal number of iterations $N$ should be set. For $n>N$ the algorithm exits without fixing a root.
\item A common problem occurs when one of the functions becomes zero before the other function. A possible way out is to substitute one of the so-fixed variables, say $x^{fin}$, in the non-zero function and to 
to run the one-dimensional M\"uller subroutine using the other variable -- $\mu(y_{n},F_2(x=x^{fin},y)) \to y^{fin}$. The algorithm then exits with a possible root: $(x^{fin},y^{fin})$.
\end{enumerate}

The procedure can be fine-tuned trough change in the starting pair of points, the initial deviation or by switching the places of the functions, or even by replacing the functions with their independent linear combinations.

The numerical experiments of the application of this algorithm on systems featuring different classes of functions can be found in \cite{arxiv}. The tests showed that the algorithm inherits some of the advantages of the one-dimensional M\"uller algorithm, like the quick convergence in proximity of the root and the vast class of functions that it can work with. The major disadvantage comes from the complicated behavior of the two-dimensional complex surfaces defined by the functions $F_{1,2}(x,y)$ which requires one to find the best combination of starting points and number of iterations in the one-dimensional M\"uller subroutine so that the algorithm converges to the required root (if it is known or suspected). Generally, it is hard to tell when one point is "close" to a root. In some cases, even if certain starting pair of points is close to a root in terms of some norm, using it as a starting point in the algorithm may still lead to convergence to another root or simply to require more iterations to reach the desired root than if other pair of starting points were used.

It is important to note that unlike Broyden's algorithm and Newton's algorithm which do not dependent on the order of the equations in the system, our two-dimensional M\"uller algorithm depends on the order of the equations. The numerical experiments show that while for some systems, changing the places of the equations has little or no effect on the convergence, in other cases, it slows down or completely breaks down the convergence. While such inherent asymmetry certainly is a weakness of the algorithm, there are ways around it. For example, one may alternate the places of the equations at each iteration or use their independent linear combinations ($F_{1,2}^{*}=\alpha_{1,2} F_1+ \beta_{1,2} F_2$). Those approaches make the algorithm more robust, but since they may cost speed, we prefer to set the order of the equations manually.

A technical disadvantage is that the whole procedure is more CPU-expensive than Newton's method and Broyden's method, since it generally makes more evaluations of the functions -- each one-dimensional M\"uller makes at least $1$ iteration on every step of the two-dimensional M\"uller, thus it makes at least $4$ evaluations of each function. This is because on each iteration of the two-dimensional M\"uller algorithm the functions in use change and thus one cannot use previous evaluations to reduce time. Still, in some cases, as demonstrated in \cite{arxiv} and also below, the so-constructed algorithm is quicker or comparable to Newton's or Broyden's method.

\section{Some applications of the method -- QNMs of nonrotating and rotating black holes}

We will work only with the confluent Heun functions, which are much better studied than the other types of Heun functions, due to their numerous physical applications. Besides their numerical implementation was used successfully in previous works by the authors. For details on the numerical testing, see \cite{arxiv}. 

The quasi-normal modes (QNMs) of a black hole (BH) are the complex frequencies ($\omega$) that govern the late-time evolution of the perturbations of the BH metric (\cite{QNM,QNM2,QNM0,QNM1,QNM21}), which have been extensively studied.

\subsection{First example: Non-rotating black hole}
First, we consider the problem of the gravitational QNMs of a nonrotating black hole so that the precision of the new method can be tested on a very well studied physical problem. The physical results in this case were published in \cite{arxiv3}, so here we will focus on the numerical details instead. 

To find the QNMs, one uses the exact solutions of the Regge-Wheeler equations, describing the linearized perturbations of Schwarzschild metric, in terms of confluent Heun functions(\cite{Fiziev1}). From \cite{Fiziev1}, when the mass of the BH is set to $2M\!=\!1$, one obtains the following system of the type \eqref{sys1}:
  \begin{align}
F_1&\!=\!(\cos(\theta) \!-\! 1) (\cos(\theta)\! +\! 1) \text{LegendreP}(l, 2, \cos(\theta))\!=\!0\label{sys} \\
F_2&\!=\text{HeunC} \Big(\!-\!2\,i\omega,2\,i\omega,4,\!-\!2\,{\omega}^{2},4\!-\!l
\!-\!{l}^{2}\!+ \!2\,{\omega}^{2},1\!-\!\mid\! r\!\mid\! e^{-i((\pi\!+\!\epsilon)/2\!+\!\it{arg}(\omega))}\! \Big)\!=\!0 \notag,\hspace{95px}
\end{align}

\noindent where $\omega$ is a complex frequency, $l$ is the angular momentum of the perturbation, $\theta \in [0,\pi]$ is the angle which we set to $\theta=\pi-10^{-7}$ and $|r|=20$. HeunC is the confluent Heun function (\cite{heun}) in \textsc{maple} notations. The parameter $\epsilon$ is a small variation ($\mid\!\epsilon\!\mid<\!1$) in the phase condition $\arg(r)\!+\!\arg(\omega)\!=\!-\!\pi/2$ (see \cite{Fiziev1}). 

Using Eqs. \eqref{sys}, we run the two-dimensional M\"uller algorithm to find the unknown $l$ and $\omega$ with precision of the algorithms set to $15$ digits.

From the theory, it is known that $l$ is an integer and $l=2, 3 ...$. Comparing with the results obtained by the two-dimensional M\"uller algorithm, for the first root $l=2$, one has $l=1.99(9)+1\times10^{-17}i$, with the first different from $9$ digit being the $17^{\text{th}}$. This shows that the new algorithm is capable of solving systems with one purely integer root in the pair with the expected precision. 

A comparison of the new algorithm with the well-known Newton's and Broyden's methods, can be found in Table \ref{table3}. 

\begin{table}[ !htb]
\vspace{-0cm}
\footnotesize
\begin{tabular}{|m{30px}|m{15px}|m{15px}|m{15px}|m{15px}|m{15px}|m{15px}|m{15px}|m{15px}|m{15px}|m{15px}|m{15px}|}
 \hline
Mode:&0&1&2&3&4&5&6&7&8&9&10\\
\hline
$t_B$ [s]&100&99&156&196&386&240&253&282&302&368&398\\
$t_{M2}$ [s]&317&413&595&741&1175&799&874&892&1364&971&1355\\
$t_{M1}$ [s]&202&218&335&357&497&457&396&613&623&594&667\\
\hline
\end{tabular}
\caption{The times needed for Broyden's method ($t_B$) and the two-dimensional M\"uller methods ($t_{M1}$ and $t_{M2}$) to fix a root. Note that while the precision of the former is $10$ digits, the precision of the other two is $14-15$ digits. To obtain those times, we solve the system: $[F_1+F_2,F_1-F_2]$ with starting points: $\omega[n]+0.01+0.01i,2.1+0.01i$, where $n=0..10$. }
\label{table3}
\end{table}

Because the phase condition  $r\!=\mid\! r\!\mid\! e^{-i(1\!+\!\epsilon)\pi/2-i arg(\omega)}$ includes the complex argument in non-analytical way, which cannot be differentiated, this problem cannot be solved directly using Newton's method. Broyden's method works but with serious limitation of its precision. This happens, because one of the roots $y$ in the pair $(x,y)$ is a real integer, while the other is complex and the algorithm fixes the integer root very quickly, thus the finite differences in the Jacobian become infinity. Because of this, the algorithm is able to fix only the first $10-11$ digits, while the other algorithms fix $14-15$ digits. Therefore, although Broyden's algorithm gives better times (see Table \ref{table3}) than the two-dimensional M\"uller algorithms, its precision is much lower and for modes with big imaginary part, it cannot be increased even by raising the software floating-point number to very high values. Furthermore, from Table \ref{table3}, one can see that the time needed for the each algorithm to exit with a root dramatically increases with $n$. This emphasize on the importance of the convergence of the algorithm, which may become critical in physical problems where multiple roots must be found (see the second example).

The numerical results for the QNMs are summed in Table (\ref{table2}). In it, the QNM frequencies obtained from Sys. \eqref{sys} are compared to those found by Andersson (\cite{Q_N_M}) with the phase amplitude method. Recently, those results were confirmed by Fiziev (see \cite{Fiziev1}) with the one-dimensional M\"uller method applied on the exact solutions of the radial equation in terms of the confluent Heun function for $l=2$ . To check the accurateness of the new method, we evaluate $\Delta=\mid \omega_{Muller2d}-\omega_{Andersson} \mid$.

\begin{table}[!h]
\vspace{-0cm}
\footnotesize
\begin{tabular}{|m{10px} | m{119px} | m{115px} | m{50px}|}
 \hline n  & Our $\omega$ & Andersson's $\omega$ & $\Delta$\\ \hline
&&&\\
0&0.7473433689+0.177924631i* &0.747343368+0.177924630i&$1.68\times 10^{-9}$\\
1&0.6934219938+0.547829750i* &0.693421994+0.547829714i&$3.60\times 10^{-8}$\\
2&0.6021069092+0.956553966i* &0.602106910+0.956553966i&$1.02\times 10^{-9}$ \\
3&0.5030099245+1.410296405i* &0.503009924+1.410296404i&$1.01\times 10^{-9}$\\
4&0.4150291596+1.893689782i* &0.415029160+1.893689782i&$4.41\times 10^{-10}$\\
5&0.3385988064+2.391216108i &0.338598806+2.391216108i&$9.67\times 10^{-10}$\\
6&0.2665046810+2.895821253i &0.266504680+2.895821252i&$1.48 \times 10^{-9}$\\
7&0.1856446684+3.407682345i &0.185644672+3.407682344i&$3.90\times 10^{-9}$\\
8&0.030649006+3.996823690i& 0+3.998000i** & 0.0306 \\
9&0.1265270180+4.605289542i&0.126527010+4.605289530i&$1.44\times 10^{-8}$\\
10&0.1531069502+5.121653272i&0.153106926+5.121653234i&$4.52\times 10^{-8}$\\
\hline
\end{tabular}
\caption{A list of the frequencies we obtained for the QNMs of Schwarzschild black hole compared with the numbers found by Andersson. $\Delta\!=\! \mid \!\omega_{Muller2d}-\omega_{Andersson} \!\mid$.  The first 5 frequencies ($n=0-4$, marked with *) were obtained also by Fiziev using the confluent Heun functions and coincide with the presented here. The 8th mode, marked with **, was obtained by Leaver \cite{Leaver}.}
\label{table2}
\end{table}

From the table, it is clear that in most cases, the modes obtained with the two-dimensional M\"uller algorithm coincide with those obtained by Andersson with more than 8 digits of precision in most cases and for modes $n=4,5,6$, there are 9 coinciding digits (Andersson published 9 digits of his frequencies). These results also confirm the roots for $n=0,1,2,3,4$ published in \cite{Fiziev1}.

The mode with biggest deviation from the expected value is $n=8$ in table \ref{table2} and it was already discussed 
in \cite{arxiv3} (and references therein). In brief its properties are due to the branch cuts in the radial function, which also lead to non-trivial dependence of the frequencies on $\epsilon$ (where $\arg(\omega)+\arg(r)=-\pi/2$ and $\epsilon<1$): for $n<4$ $\omega_n=\pm|\Re(\omega_n)|+\Im(\omega_n)i$); for $n>4$  $\omega_n(\epsilon)=-\text{sgn}(\epsilon)|\Re(\omega_n)|+\Im(\omega_n)i$ and for $n=8$, $|\epsilon|<0.75$, $\omega_{n=8}\!=\!\text{sgn}(\epsilon)\,0.030649006+3.996823690i$. 

Because of this, the value for $n=8$ in the table \ref{table2} was obtained for {\em positive} $\epsilon$ ($\epsilon=0.3$), unlike the other modes with $n\ge5$, which were obtained for $\epsilon=-0.3$. 

The frequencies presented here are stable with precision of 6 digits at the worst and usually around 9 digits with respect to a change of $\epsilon$ in the corresponding intervals. 

To confirm that the so-observed dependence of $\omega_n$ on the parameter $\epsilon$ is not a problem of the algorithm, but it is a feature of the numerical realization of the confluent Heun function, we make complex 3d plots of the radial function in use for several values of $\epsilon$.  The effect of this parameter on the branch cuts in the radial function was already discussed in \cite{arxiv3} and \cite{arxiv_Kerr}, so on Fig. \ref{fig2},  \ref{fig3}, \ref{fig4} in the Appendix we only illustrate the movement of the branch cuts under the change of $\epsilon$. 

The appearance of those branch cuts represents one more complication in the work with the confluent Heun functions, since not all of them are documented.Using the parameter $\epsilon$, however, those branch cuts can be controlled and one can obtain valuable physical results.

\subsection{Second example: Rotating black holes}

A significantly more complicated system to solve can be found in the case of QNMs from rotating black holes. In this case, one uses the exact solutions of the Teukolsky radial and angular equations, describing the linearized electromagnetic perturbations of the Kerr metric, in terms of confluent Heun functions, as stated for the first time in full detail in \cite{Fiziev1}. The two-dimensional M\"uller algorithm was applied successfully in this case too and the complete results and their discussion can be found in \cite{arxiv_Kerr}. Here, we present some details on the numerical procedures used in this case. 

From \cite{Fiziev3}, for the values of the parameters: s=-1, M=1/2, $|r|=110$, m=0, a=0.01, $\theta=\pi/3$, one obtains:
{\footnotesize
\begin{align*}
&F_{1}(x,y)=\text{HeunC}(\!-\!1.9996ix, 2.0002ix\!+\!1.0000\!, 0.0002ix\!-\!1.0000\!, 
\!-\!1.9996x(i\!+\!x), 1.9995x^2\!-\!y\!\\&\!+\!0.5000\!+\!\!1.9998ix, \!-\!110.02e^{(4.7124i\!-\!iarg(x))}\!\!+\!\!1.0000)
.{(110.00e^{(4.7124i\!-\!iarg(x))}})^{(2.00\!+\!0.0002ix)}\!=\!0\!\\
&F_{2}(x,y)= \frac{\text{HeunC}'(\!0.04x, -1.00, 1.00, -0.04x, 0.50\!-\!1.00y\!+\!0.02x\!-\!0.0001x^2, 0.25)}{\text{HeunC}(\!0.04x, -1.00, 1.00, -0.04x, 0.50\!-\!1.00y\!+\!0.02x\!-\!0.0001x^2, 0.25)}+\\ &\qquad \qquad \qquad \frac{\text{HeunC}'(\!-\!0.04x, 1.00, -1.00, 0.04x, 0.50\!-\!1.00y\!-\!0.02x\!-\!0.0001x^2, 0.75)}{\text{HeunC}(\!-\!0.04x, 1.00, -1.00, 0.04x, 0.50\!-\!y-0.02x\!-\!0.0001x^2, 0.75)}\!=\!0
 \end{align*}}
\noindent where HeunC' is the derivative of the confluent Heun function (\cite{heun}) as defined in \textsc{maple}.

For brevity, here the radial equation $F_1(x,y)$ was rounded to only 4 digits of significance. In our numerical experiments, we used the complete system with software floating-point number set to $64$, where the derivatives of the confluent Heun functions $\text{HeunC}'$ were replaced with the associate $\delta_N$ confluent Heun function according to equation (3.7) of \cite{Fiziev4}. This was done to avoid the numerical evaluation $\text{HeunC}'$ so that the peculiarities of the numerical implementation of the confluent Heun function (i.e. the use of \textsc{maple} {\em fdiff} procedure) are minimized. The difference in the times needed to fix a root when $\text{HeunC}'$ is used and when it is not used is small for the modes (i.e.$x$) with small imaginary part ($\Delta t\sim15s$), but it increases with the mode number, until it becomes significant for modes with big imaginary part (for the $10^{\text{th}}$ mode -- $R=3$ in table \ref{table0_0} -- the difference is already $\Delta t\sim 100s$). This slowdown is due to the time-consuming numerical integration and differentiation in the complex domain, needed for the evaluation of $\text{HeunC}'$. 
\begin{table}[ !h]
\centering
\vspace{0cm}
\footnotesize
\addtolength{\tabcolsep}{-1.3pt}{
 \begin{tabular}{|m{1px} |m{53px} | m{122px}  | m{37px}|m{23px} |m{23px} |}
\hline R&  $(x)^{initial}$ & $(x)^{final} $  & $t_{{}_{Broyden}}$[s] & $t_{M2}$[s]& $t_{M1}$[s]\\
  & $(y)^{initial}$ & $(y)^{final} $ & $N_{{}_{Broyden}}$ & $N_{M2}$& $N_{M1}$\\ \cline{1-6}
\hline 
\multicolumn{1}{|m{3px}|}{\multirow{3}{*}{1}}&\multicolumn{1}{c|}{$0.49+0.18i$} &$0.4965436315\!+\!0.1849695292i$&208& 102 & 92\\
&$2.001+0.1i$&$1.9999915063\!-\!0.7347653.10^{-5}i$&23&9(5)*&11(4)*\\ \cline{2-6}
\multicolumn{1}{|m{3px}|}{\multirow{3}{*}{2}}&\multicolumn{1}{c|}{$0.17+0.97i$} &$0.3495869222\!+\!1.0503235984i$&449&229&244\\
&$2.001+0.1i$&$ 2.0000392386-0.2937407.10^{-4}i$&34&12(5)*&15(5)*\\ \cline{2-6}
\multicolumn{1}{|m{3px}|}{\multirow{3}{*}{3}}&\multicolumn{1}{c|}{$0.07+5.147i$} &$0.0608496029\!+\!5.1191008697i$&868&568&489\\
&$2.001+0.051i$&$ 2.0010479243\!-\!0.2491318.10^{-4}i$&36&11(5)*&17(5)*\\
\hline
\end{tabular}}
\caption{QNMs of  Kerr BH for $s=-1$. $R$ numbers the root, $t$ and $N$ label the time and the iterations needed for the algorithms to exit. * denotes the roots dependent on the order of the equations in M1 and M2.\vspace{-0.5cm}}
\label{table0_0}
\end{table}

For that system, three pairs of starting points were used:($0.49+0.18i,2.001+0.1i$), ($0.17+0.97i,2.001+0.1i$), ($0.069+5.146i, 2.001+0.051i$). The results can be found in table \ref{table0_0}. One sees that the two modifications of the two-dimensional M\"uller algorithm M1 and M2 are much quicker than the Broyden algorithm ($t_{M1}\sim t_{M2}<t_{B}$). Newton's method once again cannot be used.

The supremacy of the M\"uller algorithms is clear and it is not isolated -- it is observed for other modes or values of the parameters (for example, for $m=1$). To check the precision of the method, the first two modes were compared with the already published results of electromagnetic QNMs of a Kerr black hole (see \cite{special31}) and were found to coincide with at least 9 digits of significance with them. We could not find a published value for the third mode. 

This example shows that in systems featuring the confluent Heun functions and their derivatives, the two-dimensional M\"uller method is much better suited than the already known algorithms.

\section{Conclusion}
The confluent Heun function appear in many physical problems. Because of the peculiarities of those functions, however, the standard numerical root-finding algorithms do not work efficiently enough on them. Here, we presented the general idea of a method for solving a system of two complex-valued nonlinear transcendental equations with complex roots based on the one-dimensional M\"uller method. This method avoids some of the problems accompanying the work with the Heun functions and it is aimed to provide adequate way to deal with systems featuring those functions. 

In the current article, the numerical results from the application of the new algorithm to systems from the QNM physics are presented. They showed that in those cases, the new method indeed works better than the standard methods. Therefore, the new method can be readily applied to find the roots of the Regge-Wheeler equation \cite{Fiziev1}, the Zerilli  equation \cite{Fiziev2}, the Teukolsky radial and angular equations \cite{Fiziev3}, all of which are solved analytically in terms of confluent Heun functions. Using this algorithm, we were able to solve {\em directly} the problem of quasi-normal modes of a Schwarzschild (\cite{arxiv3}) and Kerr black hole (\cite{arxiv_Kerr}) with higher precision than that of the Broyden method. The so found solutions agree to great extent with previous published numerical results thus confirming the usefulness of the method. 

The complete mathematical investigation of the proposed new method, and especially its theoretical order of convergence under proper conditions on the class of functions $F_1,F_2$ is still an open problem. 

For other applications of the method see the recent references \cite{PP, PP1}.

\section*{Acknowledgements}
The authors are grateful to prof. Hans Petter Langtangen for the critical reading of the early version of the manuscript and for the useful advices.

This article was supported by the Foundation "Theoretical and
Computational Physics and Astrophysics", by the Bulgarian National Scientific Fund
under contracts DO-1-872, DO-1-895, DO-02-136, and Sofia University Scientific Fund, contract 185/26.04.2010.

\section*{Author Contributions}
P.F. gave the idea and the outline of the two-dimensional generalization of the M\"uller algorithm,
chose the physical problem on which to test the algorithm and he supervised the project.

D.S. is responsible for the realization of the algorithm in
\textsc{maple} code, the testing and the optimization of the code and for the numerical results and plots presented here.

Both authors discussed the results of the tests of the algorithm, commented them and were involved
in trouble-shooting of the code at all stages. The manuscript was prepared by D.S. and edited by P.F..
\newpage
\appendix{}
\section{The epsilon-method}

\begin{figure}[!h]
 \centering
\subfigure[$\epsilon=-0.1$]{\includegraphics[scale=0.3,trim=0 0 0 8mm,clip]{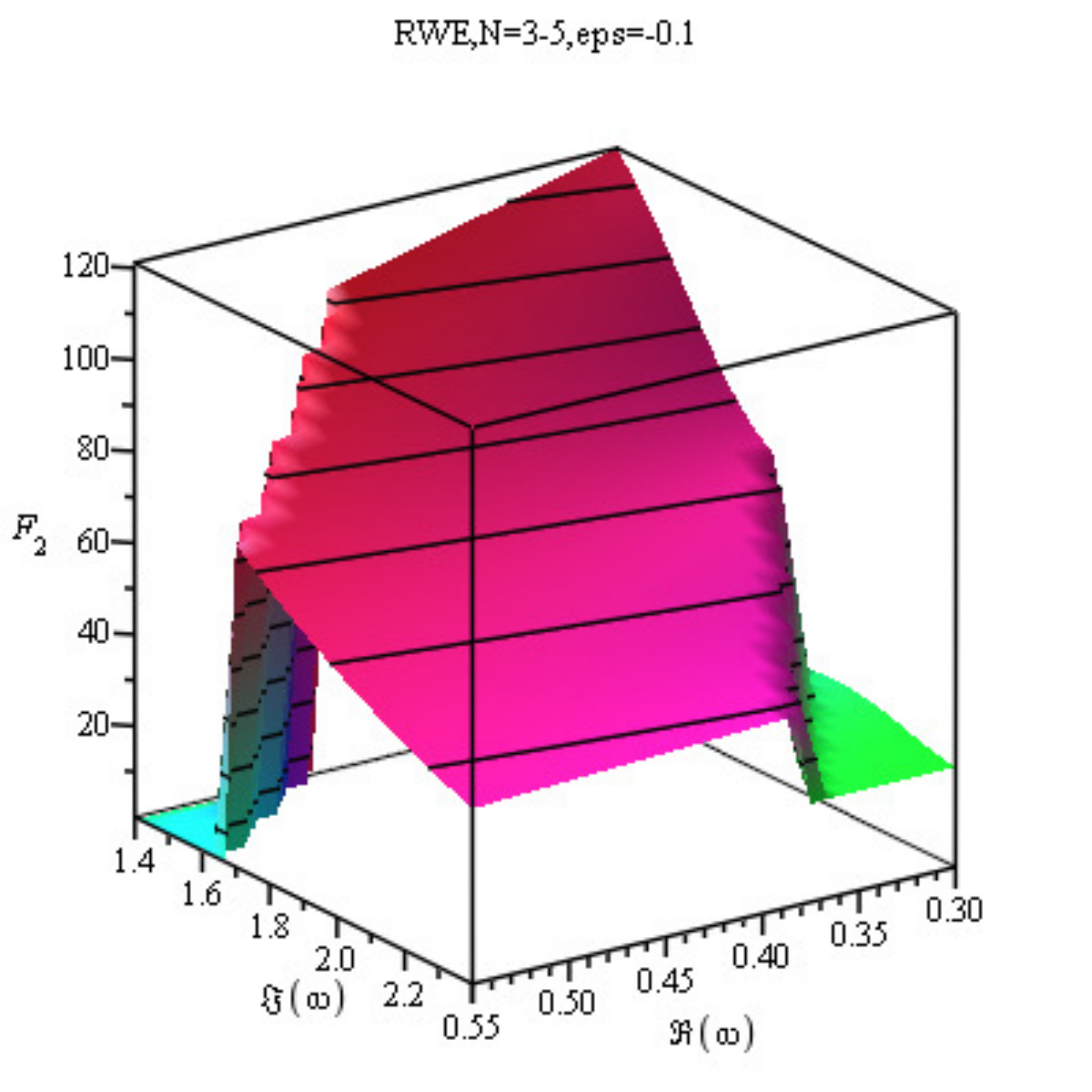}} 
\subfigure[$\epsilon=-0.05$]{\includegraphics[scale=0.3,trim=0 0 0 8mm,clip]{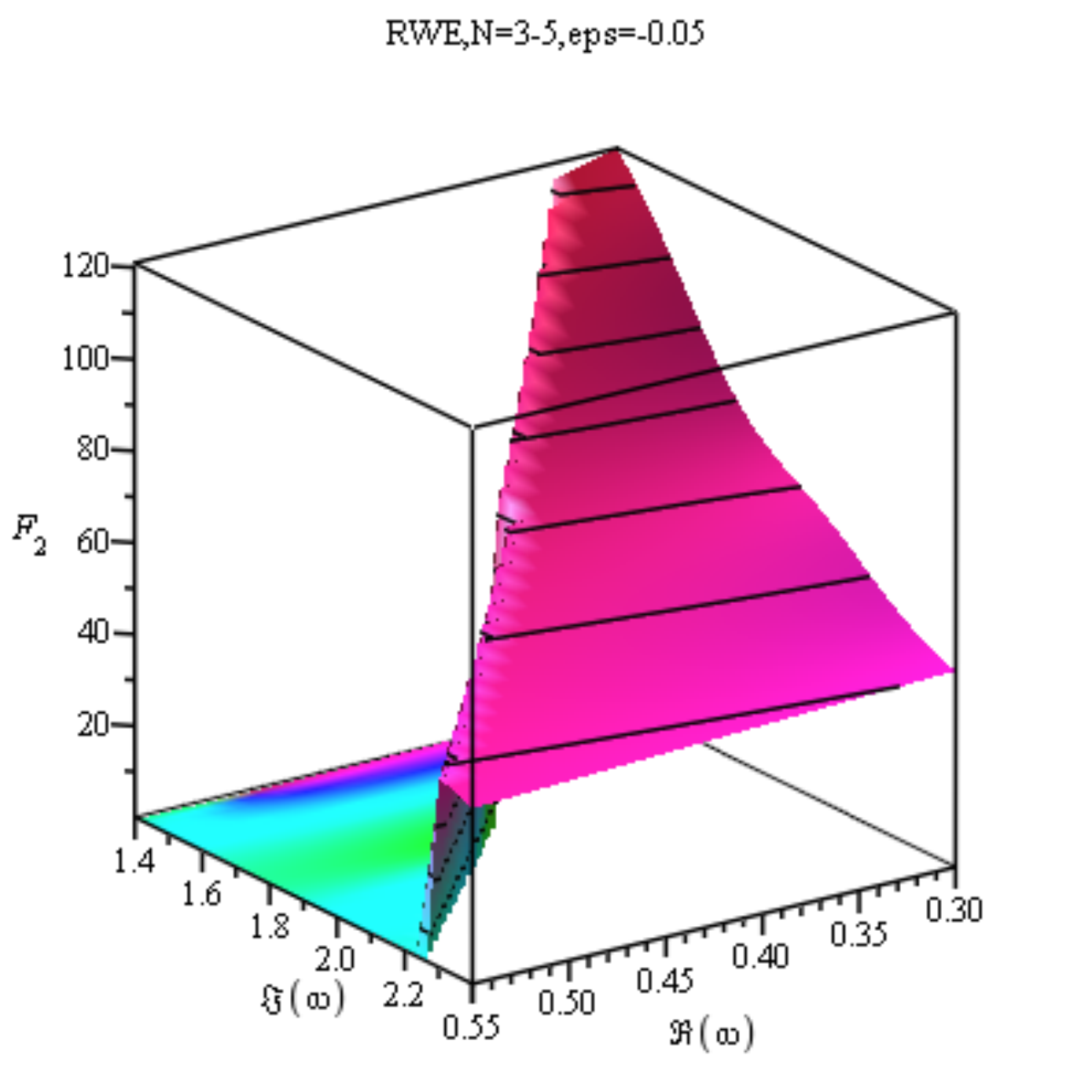}} 
\subfigure[$\epsilon=0$]{\includegraphics[scale=0.3,trim=0 0 0 8mm,clip]{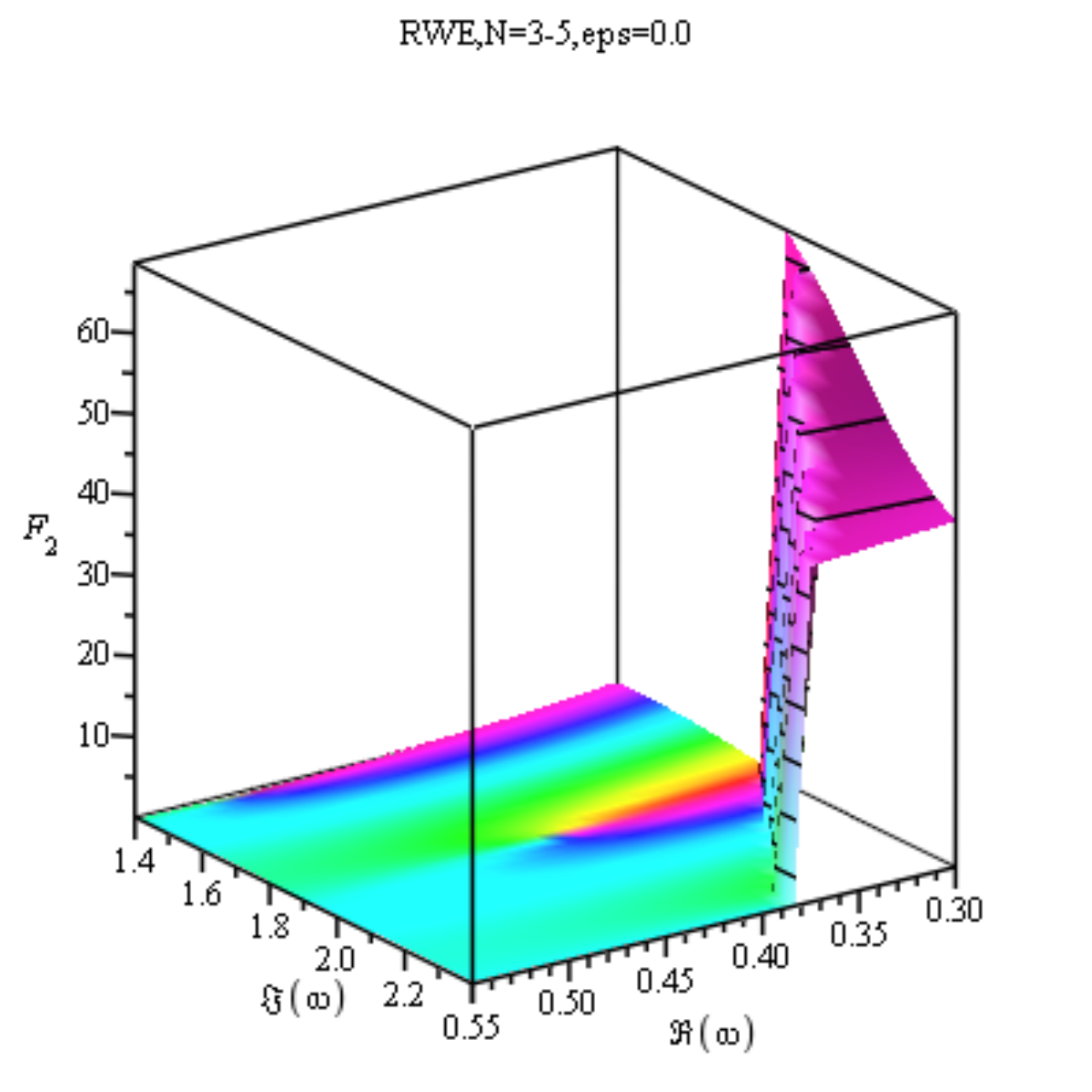}}
\subfigure[$\epsilon=0.05$]{\includegraphics[scale=0.3,trim=0 0 0 8mm,clip]{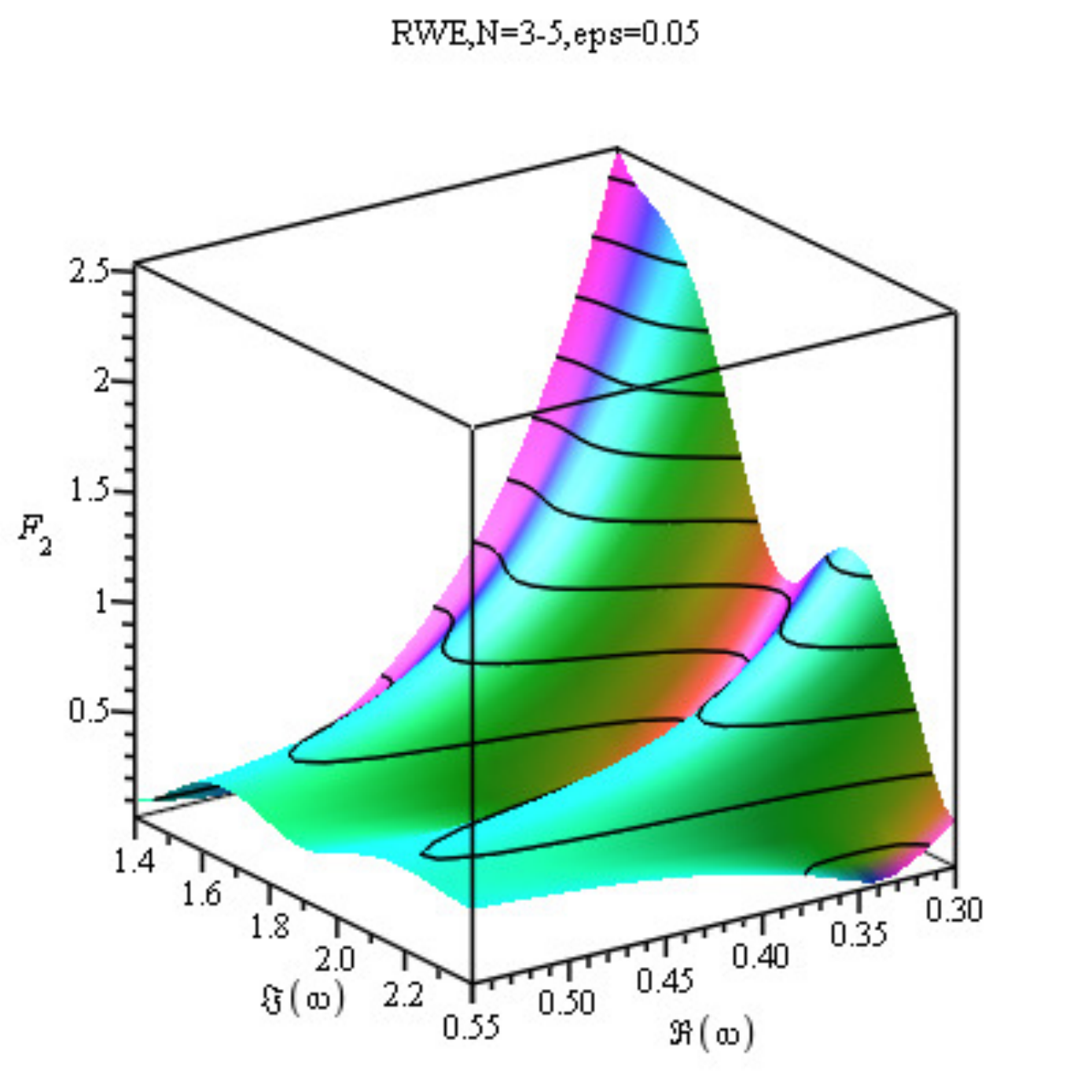}}
     \caption{ 3d plots of the function $F_2$, the solution of the RWE, in the complex interval $\omega=0.32+1.4i..0.5+2.4i$ for $\epsilon=0.05,0, -0.05,-0.1$ (the colors encode the phase of the complex function $F_2$.). The wall characteristic of the branching of the multivalued function is moved by $\epsilon$ either to the left ($\epsilon<0$) or to the right ($\epsilon>0$). Note that on Figs. \ref{fig2},\ref{fig3} and \ref{fig4}, $r=\!\mid\! r\!\mid\! e^{i((3/2\pi\!+\!\epsilon)/2\!-\!\it{arg}(\omega)}$.  }
 \label{fig2}
 \end{figure}

\begin{figure}[!ht]
\centering
\subfigure{\includegraphics[scale=0.23,trim=0 0 0 8mm,clip]{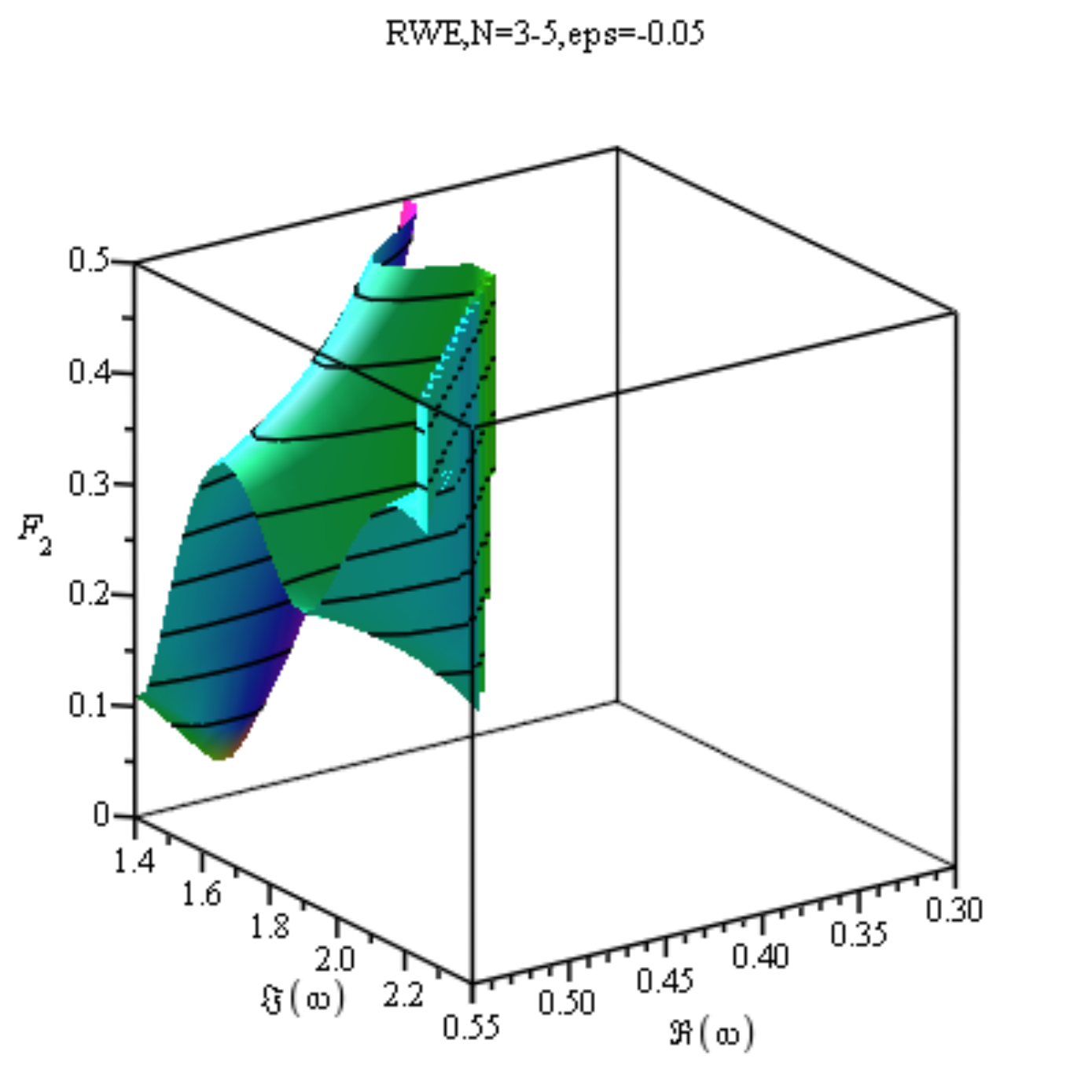}
\subfigure{\includegraphics[scale=0.23,trim=0 0 0 8mm,clip]{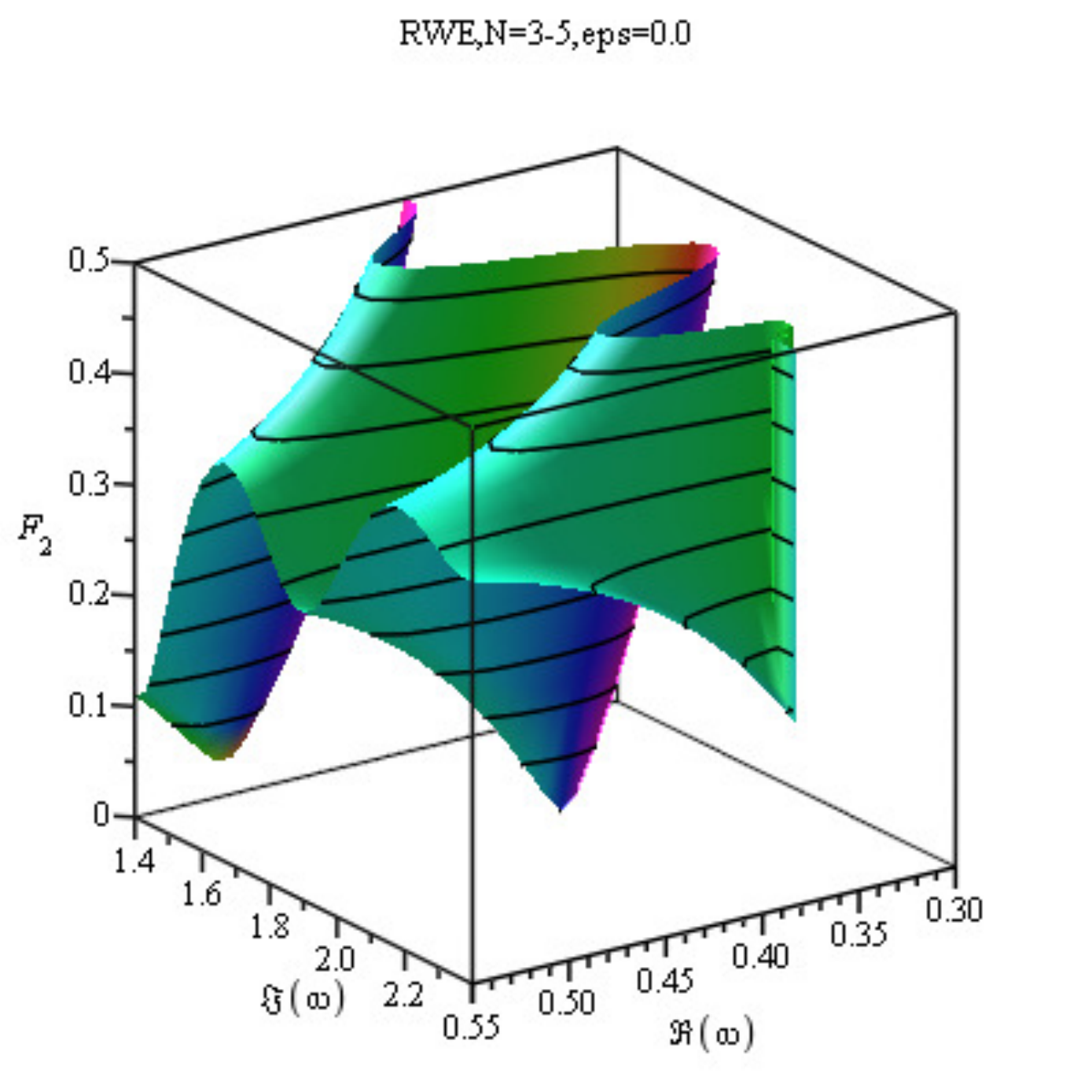}
\subfigure{\includegraphics[scale=0.23,trim=0 0 0 8mm,clip]{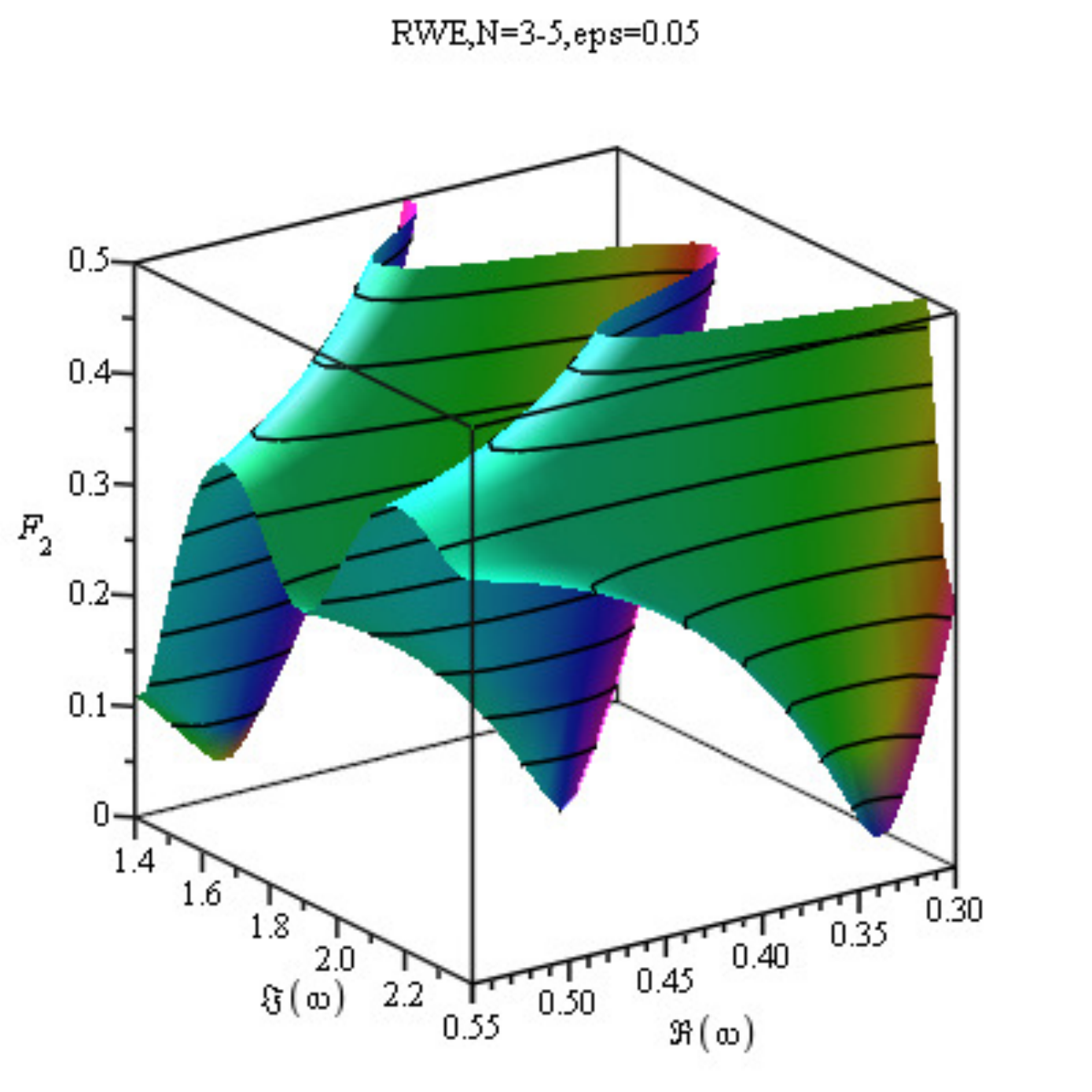}}} }
\caption{3d plots of the function $F_2$ for $\epsilon=-0.05$,$\epsilon=0$ and $\epsilon=0.05$ in the same interval for the complex $\omega$ as in \ref{fig2} scaled near the expected roots. The different $\epsilon$ lead to different profiles of the roots}
\label{fig3}
\end{figure}

\begin{figure}[!htb]
\centering
\subfigure[$\epsilon=-0.08$]{\includegraphics[scale=0.3,trim=0 0 0 8mm,clip]{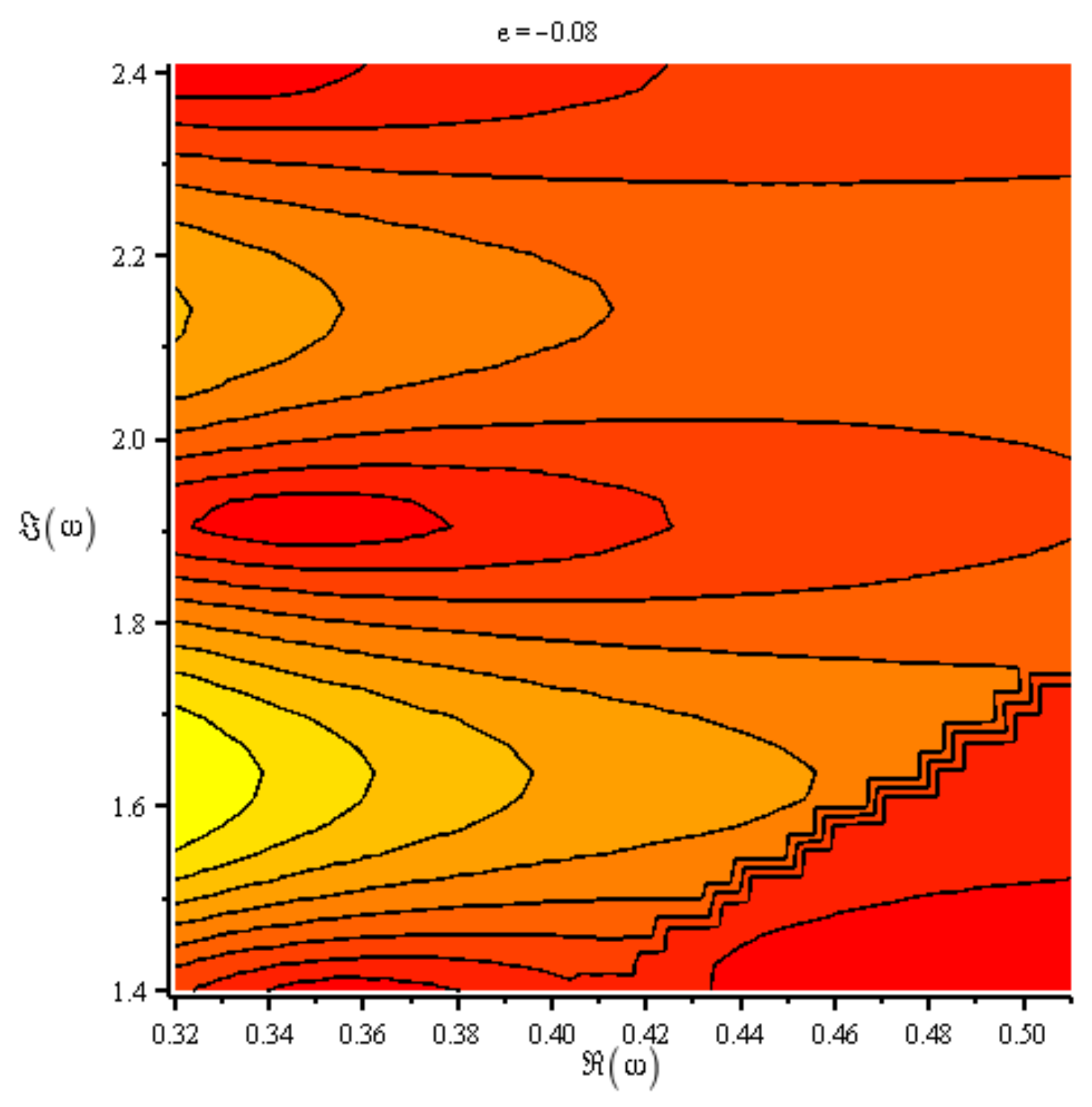}}
\subfigure[$\epsilon=-0.05$]{\includegraphics[scale=0.312,trim=0 0 0 8mm,clip]{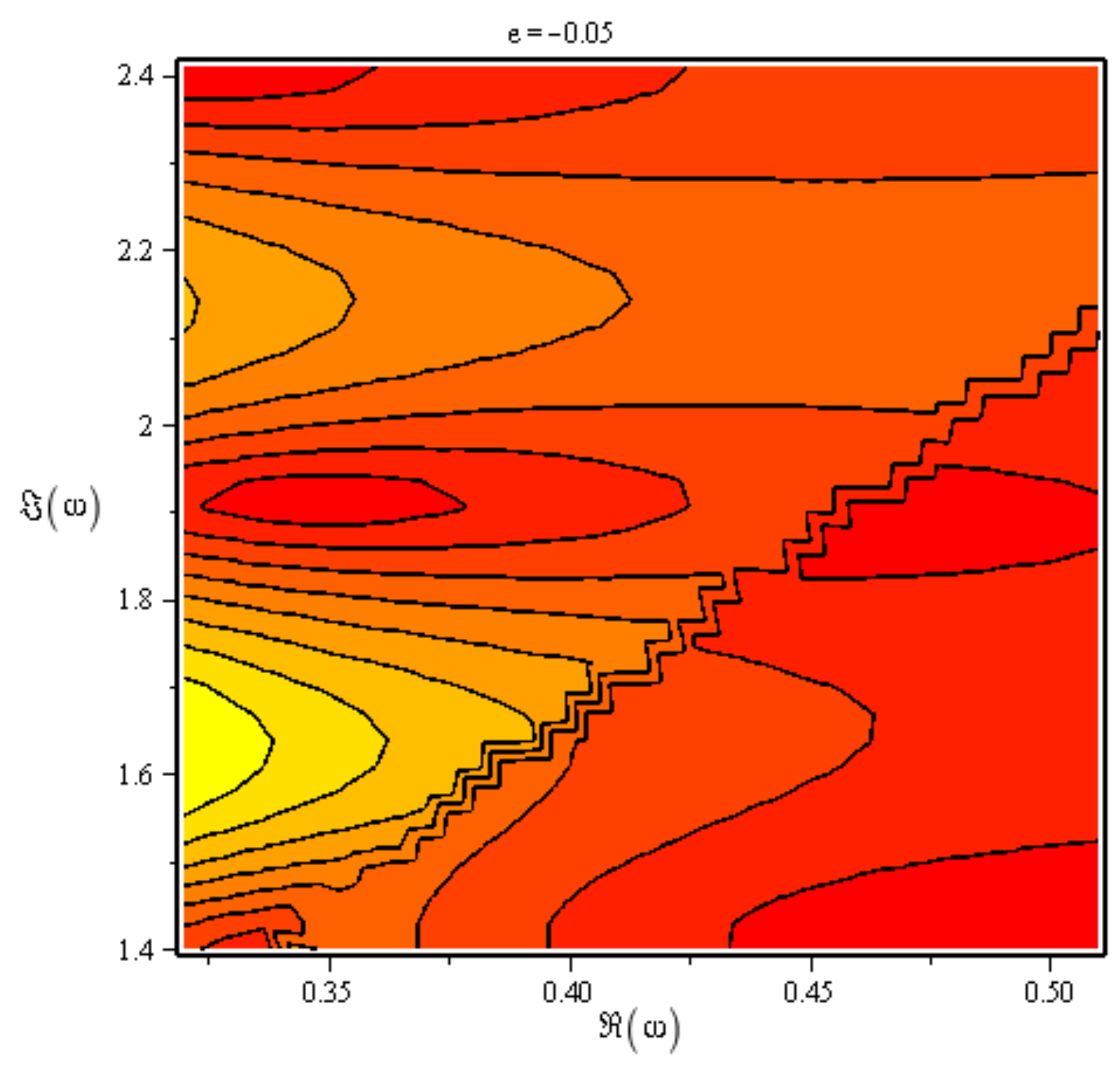}}
\subfigure[$\epsilon=0$]{\includegraphics[scale=0.3,trim=0 0 0 8mm,clip]{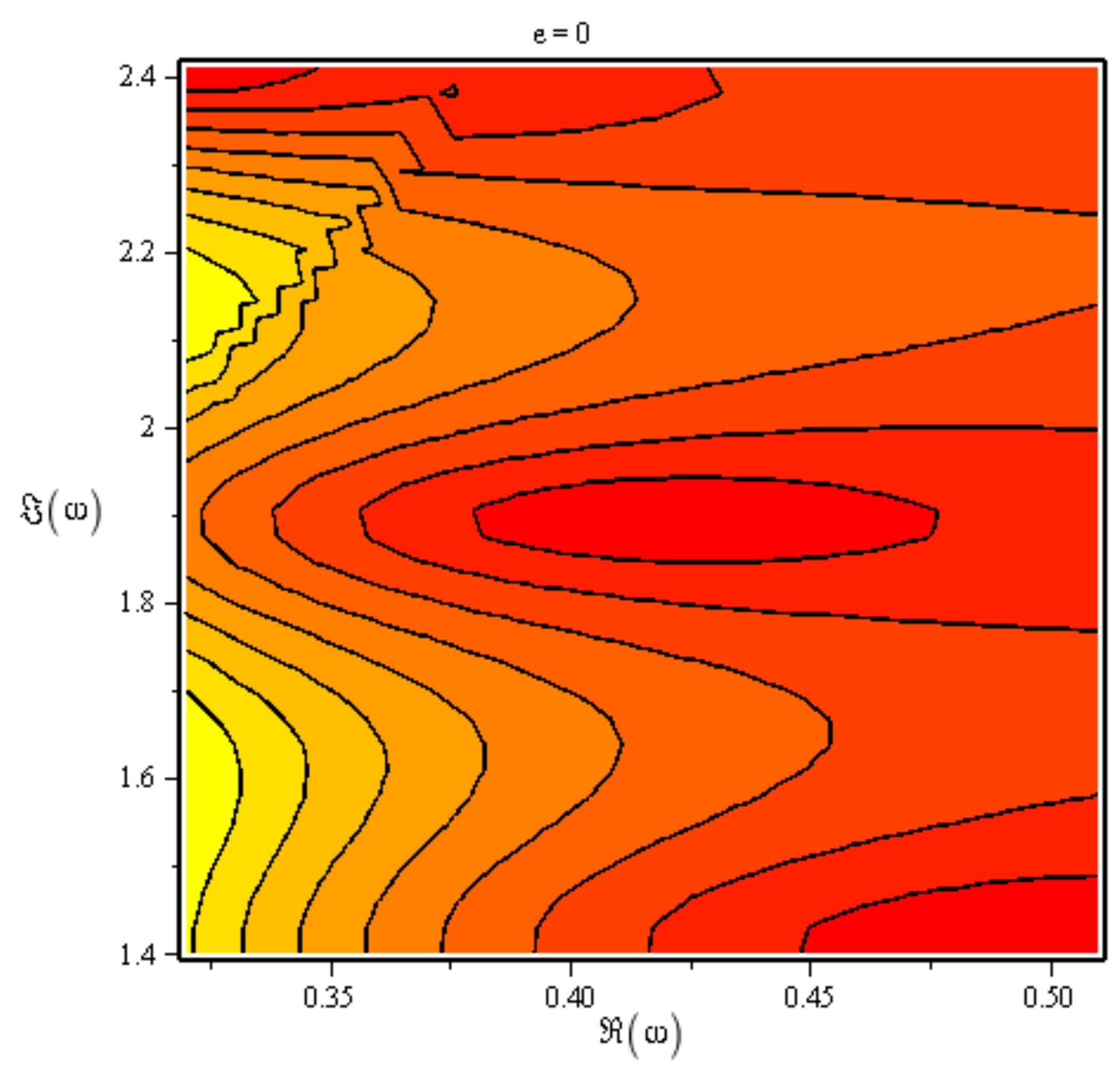}}
\subfigure[$\epsilon=0.05$]{\includegraphics[scale=0.307,trim=0 0 0 8mm,clip]{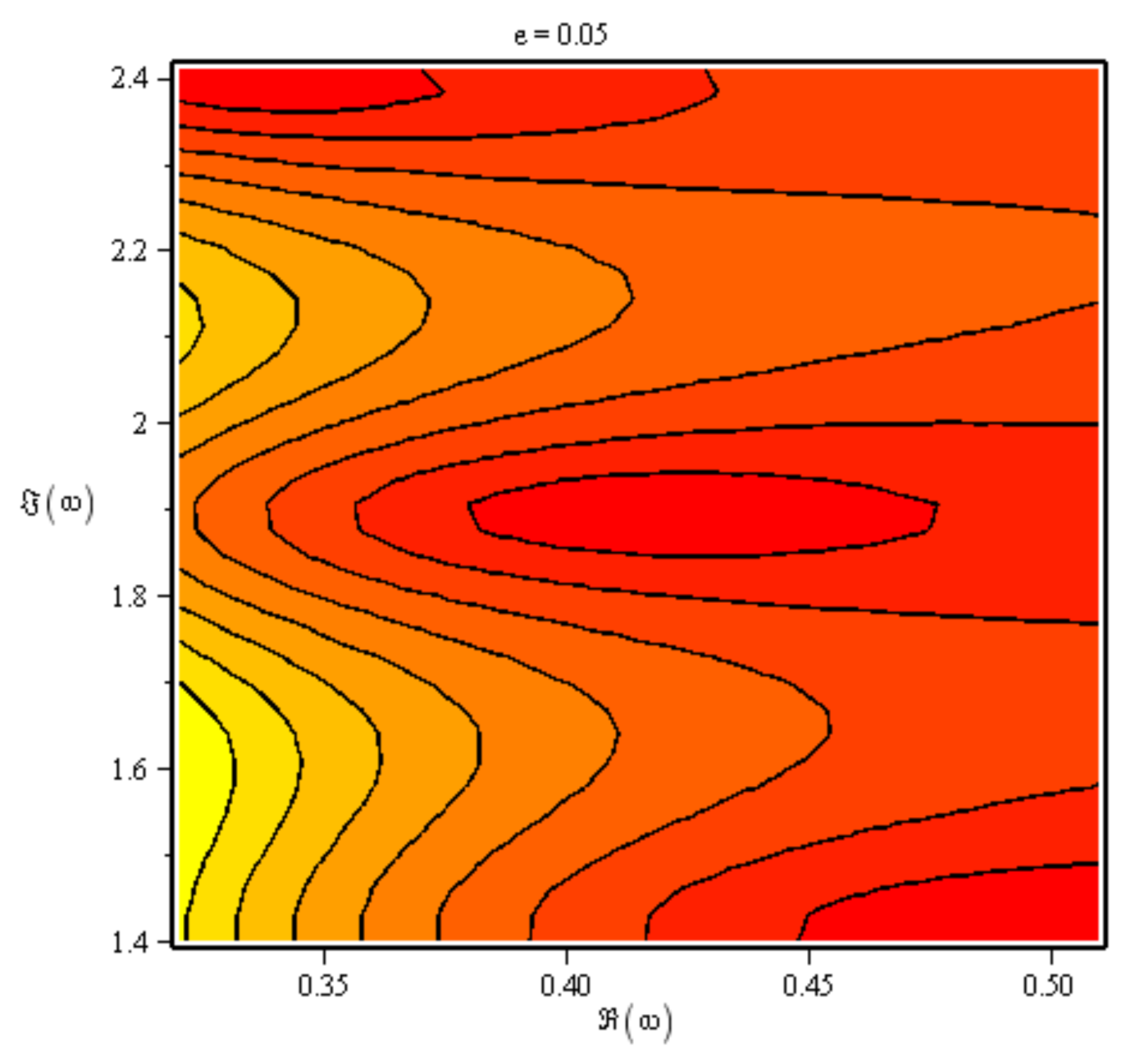}}
\caption{Plot of the level curves of the function $F_2=|R(\omega,\epsilon)|$ for $\epsilon=0.05,0,-0.05,-0.08$ for $\omega=0.32+1.4i..0.55+2.4i$. The movement of the branch cut due to the change of $\epsilon$ can be clearly seen. Note that on these plots, $F_2$ is the radial solution of the Teukolsky Radial Equation \cite{arxiv3}, but as seen from \ref{fig2}, the branch cuts for this choice of parameters coincide with those of the solution of the RWE }
\label{fig4}
\end{figure}


\newpage
\begin{thebibliography}{00}
 \bibitem{DLMF}{\em NIST Handbook of Mathematical Functions, Edited by Frank W. J. Olver, Daniel W. Lozier,
Ronald F. Boisvert, and Charles W. Clark} Cambridge University Press, 2010, 966 pages, ISBN: 978-05211-922-55, http://dlmf.nist.gov/31.17
%
\bibitem{heun3_} {\sc Slavyanov ~S.~Y., Lay ~W.}, {\em Special Functions, A Unified Theory Based on Singularities}
(Oxford: Oxford Mathematical Monographs) (2000)
%
\bibitem{heun} {\sc Heun ~K.}, Math. Ann. {\bf 33}: 161, (1889)

\bibitem{heun1_} {\sc Decarreau ~A., Dumont-Lepage ~M. ~Cl., Maroni ~P., Robert ~A. and Roneaux ~A.},  Ann. Soc. Buxelles {\bf 92}: 53, (1978)

\bibitem{heun2_} {\sc Decarreau ~A., Maroni ~P. and Robert ~A.}, 1978 Ann. Soc. Buxelles 92 151. 1995 {\em Heun's Differential Equations ed Roneaux A}, Oxford: Oxford Univ. Press, (1995)
%
\bibitem{heun1} {\sc Hortacsu ~M.}, {\em Heun Functions and their uses in Physics},   arXiv:1101.0471v1 [math-ph]
%
\bibitem{Fiziev4} {\sc Fiziev~P.~P.},  {\em Novel relations and new properties of confluent Heun's functions and their derivatives of arbitrary order}, J. Phys. A: Math. Theor. {\bf 43} (2010) 035203,    arXiv:0904.0245 [math-ph]
%
\bibitem{spectra} {\sc Staicova ~D., Fiziev ~P.}, {\em The Spectrum of Electromagnetic Jets from Kerr Black Holes and Naked Singularities in the Teukolsky Perturbation Theory}, Astrophys Space Sci {\bf 332}: 385-401 (2011), arXiv:1002.0480v2 [astro-ph.HE], (2010)
%
\bibitem{arxiv3}{\sc Fiziev ~P., Staicova ~D.},{\em Application of the confluent Heun functions for finding the QNMs of nonrotating black hole}, Phys. Rev. D {\bf 84}: 127502 (2011)  arXiv:1109.1532 [gr-qc] (2011)
%
\bibitem{arxiv_Kerr} {\sc Staicova ~D., Fiziev ~P.}, {\em New results for electromagnetic quasinormal modes of black holes}, arXiv:1112.0310v2 [astro-ph.HE]
%

\bibitem{QNM} {\sc Chandrasekhar S., and Detweiler S. L.},  {\em The quasi-normal modes of the Schwarzschild black hole},  Proc. Roy. Soc.  London A{\bf 344}: 441-452 (1975)
%
\bibitem{QNM2} {\sc Detweiler ~S.}, {\em Black holes and gravitational waves. III - The resonant frequencies of rotating holes}, ApJ:{\bf 239}: 292-295, (1980)
%
\bibitem{Muller} {\sc Press~W.~H., Teukolsky~S.~A., Vetterling~W.~T. and  Flannery~B.~P.}, {\em Numerical Recipes}, Cambridge University Press, Cambridge, England, (1992)
%
\bibitem{numerical} {\sc Forsythe ~G.~E., Malcolm ~M.~A., Moler ~C.~B.}m {\em Computer Methods for Mathematical Computations },  Prentice Hall (Prentice-Hall series in automatic computation), (1977)
%
\bibitem{Broyden} {\sc Broyden, ~C. ~G.}, {\em A Class of Methods for Solving Nonlinear Simultaneous Equations} Math. Comput. \textbf{19}, 577-593, (1965)
%
\bibitem{arxiv}
{\sc Fiziev ~P., Staicova ~D. }, {\em Two-dimensional generalization of the Muller root-finding algorithm and its applications} (2011),    arXiv:1005.5375v2 [cs.NA]
%
\bibitem{Muller2} {\sc M\"uller, David ~E.}, {\em A Method for Solving Algebraic Equations Using an Automatic Computer}, MTAC {\bf 10}:208-215 (1956)
%
\bibitem{Fiziev1} {\sc Fiziev P.~P.}, {\em Exact Solutions of Regge-Wheeler Equation and Quasi-Normal Modes of Compact Objects}, Class. Quant. Grav. {\bf 23} 2447-2468 (2006), arXiv:0509123v5 [gr-qc]
%

\bibitem{Q_N_M} {\sc Andersson~N.}, {\em A numerically accurate investigation of black-hole normal modes},
Proc. Roy. Soc. London A\textbf{439} no.1905: 47-58 (1992)
%
\bibitem{special31} {\sc Berti ~E., Cardoso ~V., Will ~C.~M. }, {\em On gravitational-wave spectroscopy of massive black holes with the space interferometer LISA}, Phys.Rev.D {\bf 73}:064030, (2006), arXiv:0512160v2 [gr-qc]
%
\bibitem{QNM0}
    {\sc Chandrasekhar ~S.}, {\em The mathematical theory of black holes}, , Clarendon Press/Oxford University Press (International Series of Monographs on Physics. Volume 69), (1983)
%
\bibitem{QNM1}
{\sc Ferrari~V., Gualtieri L.}, {\em Quasi-normal modes and gravitational wave astronomy}, Gen.Rel.Grav.{\bf 40}: 945-970 (2008), arXiv:0709.0657v2 [gr-qc] (2007)
%
\bibitem{QNM21}
{\sc Konoplya, ~R. ~A., Zhidenko, ~A.}, {\em Quasinormal modes of black holes: from astrophysics to string theory}, to be published in Reviews of Modern Physics, arXiv:1102.4014v1 [gr-qc] (2011)
%
\bibitem{Fiziev3} {\sc Fiziev~P.~P.},  {\em Classes of exact solutions to the Teukolsky master equation}
                       Class. Quantum Grav. {\bf 27}  135001 (2010), arXiv:0908.4234v4  [gr-qc]
%
\bibitem{Leaver} {\sc Leaver ~E. ~W.}, {\em An analytic representation for the quasi-normal modes of Kerr black holes}, Proc. Roy. Soc. London A{\bf 402}: 285-298 (1985)
%
\bibitem{special3} {\sc Berti~E., Cardoso~V. and Starinets~A.~O.}, {\em Quasinormal modes of black holes and black branes}, Class. Quantum Grav. {\bf 26} 163001 (108pp) (2000)
%
\bibitem{special1} {\sc Berti ~E.}, {\em Black hole quasinormal modes: hints of quantum gravity?},     arXiv:gr-qc/0411025v1 (2004)
%
\bibitem{Fiziev2} {\sc Fiziev~P.~P.},  {\em Teukolsky-Starobinsky identities: A novel derivation and generalizations}, Phys. Rev. D{\bf 80}, 124001 (2009), arXiv:0906.5108 [gr-qc]
%
\bibitem{PP}
{\sc Cadoni ~M., Pani ~P.}, {\em Holography of charged dilatonic black branes at finite temperature}, 	JHEP {\bf 1104}:049,2011, arXiv:1102.3820v2 [hep-th]] (2011)

\bibitem{PP1}
{\sc Pani ~P.}, {\em Applications of perturbation theory in black hole physics. [Doctoral Thesis]} (2011), http://veprints.unica.it/553/

\end{thebibliography}
\end{document}